\documentclass[11pt]{article}
\usepackage{enumerate}
\usepackage{romannum}

\usepackage{enumitem}
\usepackage{fixltx2e}
\usepackage{graphicx}
\usepackage{longtable}
\usepackage{float}
\usepackage{setspace}
\usepackage{wrapfig}
\usepackage{soul}
\usepackage{amsmath}
\RequirePackage{amsthm}
\usepackage{mathtools}
\usepackage{mathrsfs}
\usepackage{textcomp}
\usepackage{marvosym}
\usepackage{wasysym}
\usepackage{latexsym}
\usepackage{amssymb}
\usepackage{hyperref}
\usepackage{graphicx}
\usepackage{subfig}
\usepackage{longtable}
\usepackage{float}
\usepackage{array}
\usepackage[ruled]{algorithm}
\usepackage[noend]{algpseudocode}
\usepackage[table]{xcolor}
\usepackage{arydshln}
\usepackage{cancel}
\usepackage{url}
\usepackage{comment}
\usepackage[square, comma]{natbib}
\setcitestyle{numbers}
\usepackage[margin=1in]{geometry}
\usepackage[affil-it]{authblk}
\usepackage{titlesec}
\usepackage[flushleft]{threeparttable}
\usepackage{color}
\titleformat{\paragraph}
{\normalfont\normalsize\bfseries}{\theparagraph}{1em}{}
\titlespacing*{\paragraph}
{0pt}{3.25ex plus 1ex minus .2ex}{1.5ex plus .2ex}

\newcommand{\subscript}[2]{$#1 _ #2$}

\allowdisplaybreaks

\numberwithin{equation}{section}

\newtheorem{theorem}{Theorem}[section]

\newtheorem{lemma}{Lemma}[section]
\newtheorem{definition}[theorem]{Definition}
\newtheorem{corollary}[theorem]{Corollary}

\input{./Definitions}

\allowdisplaybreaks

\newcommand{\fY}{\mathcal{Y}}

\title{Asymptotic Normality of the Posterior Distributions in a Class of Hidden Markov Models}

\author[1]{Chunlei Wang \thanks{\url{chunlei-wang@uiowa.edu}}}
\author[1]{Sanvesh Srivastava \thanks{\url{sanvesh-srivastava@uiowa.edu}}}

\affil[1]{Department of Statistics and Actuarial Science, The University of Iowa}

\date{\today}

\begin{document}
\pagenumbering{arabic}
\maketitle

\begin{abstract}
We show that the posterior distribution of parameters in a hidden Markov model with parametric emission distributions and discrete and known state space is asymptotically normal. The main novelty of our proof is that it is based on  a testing condition and the sequence of test functions is obtained using an optimal transportation inequality. 
\end{abstract}
  
\noindent%

\section{Introduction}

Hidden Markov models with parametric emission distributions and discrete and known state space, which we shorten as HMMs, have been widely used in diverse areas for analyzing dependent discrete-time data \citep{fru06, elliott2008hidden,zucchini2017hidden}. Our focus is on studying asymptotic properties of the posterior distribution of parameters in an HMM \citep{cappe2006inference,de2008asymptotic}. We prove a Bernstein-von Mises theorem for HMMs using a testing condition.  Informally, this theorem says that the posterior distribution of parameters in an HMM converges in distribution to a normal distribution that centers at the maximum likelihood estimator of parameters, and its covariance matrix equals inverse of the Fisher information matrix.

Consider the setup of a general hidden Markov model as a discrete-time stochastic process $(X_{t},Y_{t})_{t \geq 1}$, where $(Y_{t})_{t\geq 1}$ process is observed and $(X_{t})_{t\geq 1}$ process is unobserved and Markov. Its state space includes all possible values of $(X_{t})_{t\geq 1}$ and emission distribution is defined as the conditional distribution of $Y_{t}$ given $X_{t}$. Our HMMs are a special case of this setup in that the emission distributions are parametric and state space is discrete and known. There are many theoretical results about frequentist estimation in hidden Markov models with discrete, continuous, or unknown state spaces and parametric or non-parametric emission distributions \citep{baum1966statistical, leroux1992maximum,bickel1998asymptotic,jensen1999asymptotic, douc2001asymptotics,douc2011consistency}. The equivalent asymptotic Bayesian results have started to emerge recently \citep{de2008asymptotic,gassiat2014posterior,vernet2015posterior,gassiat2018efficient,douc2020posterior}.

Our theoretical results are closest to those of \citet{de2008asymptotic}. They prove a Bernstein-von Mises theorem for HMMs using Taylor expansion of the log likelihood function. Our proof differs from theirs in that we adopt the setup of \citet{ghosal2000convergence} for deriving posterior contraction rates in general Bayesian models. Using this setup, we adapt the proof of Theorem 10.1 in \cite{Van00} to show that the posterior distribution of parameters in an HMM is asymptotically normal. Compared to \cite{de2008asymptotic}, our assumption on the log likelihood function is milder in that our proof assumes only the first order smoothness of the log likelihood function. More importantly, our proof relies on the local asymptotic normality condition for HMMs, which is satisfied under the assumptions of \citet{bickel1998asymptotic,de2008asymptotic}.
Our proof techniques are similar to those used for investigating asymptotic properties of the posterior distribution in extensions of HMMs with  non-parametric emission distributions \citep{gassiat2014posterior,vernet2015posterior}. A key step in our proof is to establish a testing condition based on a sequence of test functions. We are not aware of any existing test function that is tuned for such applications in HMMs.

Our main contribution is to establish the testing condition using a sequence of test functions obtained from an optimal transportation cost inequality. The testing condition measures the complexity of the HMM family and establishing it is a key step in showing the consistency of posterior distribution \citep{ghosal2000convergence}. Our construction of the sequence of test functions is based on the $L^1$-transportation cost information inequality for HMMs \citep{djellout2004transportation, kontorovich2008concentration} and Lipschitz concentration of log likelihood function of HMMs \citep{le2000exponential}. \citet{hu2011transportation} has constructed the likelihood ratio test with exponentially decaying error probabilities for testing simple hypotheses. Extending this work, we construct the sequence of test functions with exponentially decaying
error probabilities for testing general hypotheses using the $L^1$-transportation cost information inequality. The sequence of test functions satisfies the testing condition in \cite{ghosal2000convergence,Van00} and is used to prove the Bernstein-von Mises theorem for HMMs.

This paper is structured in four sections. In Section \ref{sec:pre}, we state the HMM setup, define the prediction filter, and restructure the results of \citet{kontorovich2008concentration} and \citet{hu2011transportation} 
in the context of this paper. The assumptions and main results about the testing condition are presented in Section \ref{sec:main}. The contraction rate and asymptotic normality of the posterior distribution of parameters in an HMM follow immediately from the testing condition. Finally, Section \ref{sec:proof} contains proofs of the main results.

\section{Setup and Background}
\label{sec:pre}

Consider the HMM setup in greater detail. Recall that the HMM is a discrete-time stochastic process $(X_{t},Y_{t})_{t \geq 1}$ and that $(Y_{t})_{t\geq 1}$ process is observed and $(X_{t})_{t\geq 1}$ process is unobserved and Markov. 
The hidden Markov chain $(X_{t})_{t \geq 1}$ has state space $\mathcal{X} = \{1,\ldots,S\}$ with a known $S \in \mathbb{N}$. The $(X_{t})_{t \geq 1}$ process has an initial stationary distribution $r(a) = \mathbb{P}(X_{1} = a)$ for $a \in \mathcal{X}$ and is time homogeneous with transition kernel
\begin{align}
\label{eq:1}
Q_{ab} = q(a,b) = \mathbb{P}(X_{t+1} = b | X_{t} = a),\quad t \geq 1,\;a,b\in\mathcal{X}.
\end{align}
Given $(X_{t})_{t\geq 1}$, $(Y_{t})_{t \geq 1}$ is a sequence of independent random variables on the metric space $(\mathcal{Y},d_{\mathcal{Y}})$ and the conditional distribution of $Y_{t}$ depends only on $X_{t}$. Furthermore, the conditional distribution of $Y_{t}$ given $X_{t} = x_{t}$ does not depend on $t$ and has a density function $g(\cdot \mid x_{t})$ with respect to the $\sigma$-finite measure $\mu$ on $\mathcal{Y}$. Let $\Theta\subseteq \mathbb{R}^{d}$ denote the paramter space with fixed dimension $d \in \mathbb{N}$. The density functions $r(\cdot)$, $q(\cdot,\cdot)$, and $g(\cdot\mid \cdot)$ belong to a parametric family indexed by $\theta$, which is denoted as $\{r_\theta(\cdot), q_\theta(\cdot,\cdot), g_\theta(\cdot\mid \cdot), \theta \in \Theta \subseteq \mathbb{R}^{d}\}$. Our focus is Bayesian inference on $\theta$ given the observed data. 

Consider the joint and marginal distributions of the augmented and observed data. Let $n$ be the number of observations, $(Y_{1},\ldots,Y_{n})$ be the observed data, $(X_{1},\ldots,X_{n},Y_{1},\ldots,Y_{n})$ be the augmented data, and $c_{i}^j$ be the shorthand for sequence $c_i, \ldots, c_j$ for $i \leq j$. Then, the joint and marginal densities of the augmented and observed data under
$(\text{counting measure})^{n}\otimes \mu^{n}$ and $\mu^{n}$, respectively, are
\begin{align}
p_{\theta}(x_{1}^n, y_{1}^n) &= r_{\theta}(x_{1})\prod_{k=1}^{n-1}q_{\theta}(x_{k},x_{k+1})g_{\theta}(y_{k}|x_{k}) g_{\theta}(y_{n}|x_{n}), \nonumber\\
p_{\theta}(y_{1}^n) &= \sum_{(x_{1}^{n})\in \mathcal{X}^{n}} p_{\theta}(x_{1}^n,y_{1}^n),
\label{eq:lik:f}
\end{align}
where $x_{1}^n \in \mathcal{X}^n$, $y_{1}^n \in \mathcal{Y}^n$, and $\theta \in \Theta$. Denote the distribution of $Y_1^n$ as $\mathbb{P}_{\theta}^{(n)}$, density $p_{\theta}(y_{1}^n)$ as $p_{\theta}^{(n)}$, and the expectation under $\PP^{(n)}_{\theta}$ as $\mathbb{E}_{\theta}^{(n)}$. If $\theta_0 \in \Theta$ is the true parameter, then $\PP^{(n)}_{\theta_0}$ is the true distribution of $Y_{1}^{n}$.

We now define related conditional distributions and the \emph{extended} HMM sequence. For $2 \leq t \leq n$ and $\theta \in \Theta$, the conditional distribution of $X_{t}$ given $Y_{1}^{t-1}$ is called the \emph{prediction filter} at time $t$ and is defined as
\begin{align}
\label{eq:filter}
p_{t}^{\theta}\vcentcolon = \left\{ \mathbb{P}_{\theta}(X_{t}= 1 \mid Y_{1}^{t-1}), \ldots, \mathbb{P}_{\theta}(X_{t}= S \mid Y_{1}^{t-1}) \right\}^\mathsf{T} \in \mathcal{E} ,
\end{align}
where $\mathcal{E}$ is the space of all probability distributions on $\mathcal{X} $ equipped with total variation distance $\|\cdot\|_{\text{TV}}$. Baum's forward equation further implies that for $t = 1, \ldots, n-1$,
\begin{align}
\label{eq:fb}
p_{t+1}^{\theta} (\cdot) &= \frac{\PP_{\theta}(X_{t+1} = \cdot, Y_{t}\mid Y_{1}^{t-1})}{\PP_{\theta}(Y_{t}\mid Y_{1}^{t-1})}
= \frac{\sum_{x_{t}}q_{\theta}(x_{t}, \cdot) g_{\theta}(Y_{t} \mid X_{t} = x_{t}) p_{t}^{\theta}(x_{t}) }{\sum_{x_{t}}g_{\theta}(Y_{t} \mid X_{t} = x_{t}) p_{t}^{\theta}(x_{t})} \nonumber\\
&\vcentcolon = f^{\theta}(Y_{t},p_{t}^{\theta}),
\end{align}
where $f^{\theta}$ is a measurable function on $\mathcal{Z} = \mathcal{Y} \times \mathcal{E}$ equipped with metric $d_{\mathcal{Z}} = d_{\mathcal{Y}} + \|\cdot\|_{\text{TV}}$. Extending the definition of $f^{\theta}(Y_{t},p_{t}^{\theta})$ for lags $s = 0, 1, \ldots, t - 1$, we recursively define $p_{t+1}^{\theta}$ for $s = 0,
\ldots, t-1$ as
\begin{align}
\label{eq:funcf}
p_{t+1}^{\theta} = f_0^\theta(Y_t^t, p^{\theta}_{t}) = f_{1}^{\theta}(Y_{t-1}^{t},p^{\theta}_{t-1}) = \cdots = f_{s}^{\theta}(Y_{t-s}^{t},p^{\theta}_{t-s}) = \cdots = f_{t-1}^{\theta}(Y_{1}^{t},p^{\theta}_{1}),
\end{align}
where $f_{0}^{\theta} =f^{\theta}$ and $p_{1}^{\theta} = r_{\theta}$.
The \emph{extended HMM} sequence is defined as $(Y_{t},p^{\theta}_{t})_{t = 1}^{n} \in (\mathcal{Y}\times \mathcal{E})^{n} = \mathcal{Z}^{n}$, where $(p^{\theta}_{t})_{t=1}^{n}$ is defined in \eqref{eq:funcf}. Let $l_{n}(\theta,Y_{1}^{n}) = \log p_{\theta}(Y_{1}^n)$ be the log likelihood function of $\theta$, where $p_{\theta}(y_{1}^n)$ is defined in \eqref{eq:lik:f}. Using the extended HMM sequence, $l_{n}(\theta,Y_{1}^{n})$ is expressed  as
\begin{align}
l_{n}(\theta,Y_{1}^{n}) &= \log p_{\theta}(Y_{1})+\log \prod_{t=2}^{n}p_{\theta}(Y_{t}|Y_{t-1},\ldots,Y_{1})
\nonumber\\
& =
\log \left\{\sum_{x_1}g_{\theta}(Y_{1}|x_{1})r_{\theta}(x_1)\right\} +\sum_{t=2}^{n}\log\left\{ \sum_{x_t}g_{\theta}(Y_{t}|x_{t})\PP_{\theta}(X_{t} = x_{t}| Y_{t-1},\ldots,Y_{1})\right\}
\nonumber\\
& = \sum_{t=1}^{n} \log \sum_{x_{t}}p_{t}^{\theta}(x_{t})g_{\theta}(Y_{t}|x_{t}).
\label{eq:lik:a}
\end{align}
The log likelihood function in \eqref{eq:lik:a} depends on $(Y_{t},p^{\theta}_{t})_{t = 1}^{n}$. The main advantage of
\eqref{eq:lik:a} is that the log likelihood function $l_{n}(\theta,Y_{1}^{n})$ can be expressed in an additive form and easier to manipulate. Unfortunately, this is not possible if we use $p_{\theta}(y_{1}^n)$ in \eqref{eq:lik:f} to define  $l_{n}(\theta,Y_{1}^{n})$.

We now state the optimal transportation inequality and related results that are required for constructing the sequence of test functions.
\begin{definition}
\label{def:t1}
Let $(E,d)$ be a metric space. A probability measure $\mu$ on $(E,d)$ satisfies $L^{1}$-transportation cost information inequality if for any probability measure $\nu$ on $(E,d)$, there exists a constant $C>0$ such that
\begin{align}
\label{eq:t1}
W_{1}^{d}(\mu,\nu) \leq \sqrt{2CH(\nu\mid \mu)},
\end{align}
where $W_{1}^{d}$ is the optimal transportation cost with cost function $d(\cdot,\cdot)$ and $H(\cdot\mid \cdot)$ denotes the Kullback-Leibler divergence. If a measure $\mu$ satisfies \eqref{eq:t1}, we denote this relation as $\mu \in T_{1}(C)$.
\end{definition}
\noindent The $L^{1}$-transportation cost information inequality is important because it is an equivalent characterization of sub-Gaussian measure concentration.
\begin{theorem}\cite[Theorem 3.1]{bobkov1999exponential}
\label{thm:ti}
A measure $\mu$ defined on a metric space $(E,d)$ satisfies $T_{1}(C)$-inequality if and only if for any $\mu$-integrable Lipschitz function $F: E\rightarrow \mathbb{R}$ and any $\lambda \in \mathbb{R}$,
$$\int_{E} e^{\lambda(F - \langle F\rangle_{\mu})} d\mu \leq e^{\tfrac{\lambda^{2}}{2}C\|F\|^{2}_{\text{\text{Lip}}}},$$
where $\langle F\rangle_{\mu} = \int_{E}F d\mu$ and $\|F\|_{\text{\text{Lip}}}$ is the Lipschitz constant of $F$. In this case, for any $r > 0$,
\begin{align}
\label{eq:conc-ineq}
\mu(F - \langle F\rangle_{\mu}>r) \leq e^{-\frac{r^{2}}{2C\|F\|^{2}_{\text{Lip}}}}.
\end{align}
\end{theorem}

An important property of $L^{1}$-transportation cost information inequality is called the \emph{tensorization principle}. For $i=1,\dots,n$, let $\mu_{i}$ be a probability measure on metric space $(E_{i},d_{i})$ such that $\mu_{i} \in T_{1}(C)$. Consider the product measure $\mu_1\otimes \cdots\otimes \mu_{n}$ on $(E_{1}\times \cdots \times E_{n}, d_{l_1})$, where the $l_{1}$-metric is defined as
\begin{align}
\label{eq:metric1}
d_{l_{1}}(u,v) = \sum_{i=1}^{n} d_{i}(u_{i},v_{i}), \;
u = (u_{1},\ldots,u_{n}),\; v = (v_{1},\ldots,v_{n}),\; u_i, v_i \in E_i.
\end{align}
Then, the tensorization principle for $L^{1}$-transportation cost information inequality states that
\begin{align}
\label{eq:tensor}
\mu_1 \otimes \cdots \otimes \mu_{n} \in T_{1}(nC);
\end{align}
see \cite[Chapter 2]{van2014probability} for greater details. For a product measure, this property is called \emph{independent tensorization}.
Under mild conditions, the probability measure $\mathbb{P}^n_\theta$ induced by $(Y_{t})_{t =1}^{n}$ satisfies the $T_{1}(nC_{H}^{\theta})$-inequality on the metric space $(\mathcal{Y}^{n},d_{l_{1}})$ for some constant $C_{H}^{\theta}$. If $\mathcal{Y}$ is countable, then the constant $C_{H}^{\theta}$ is determined by the mixing property of the hidden Markov chain $(X_{t})_{t \geq 1}$ \citep{kontorovich2008concentration}. Furthermore, $C_{H}^{\theta}$ is also related to the measure concentration of emission distribution if $\mathcal{Y}$ is uncountable \citep{hu2011transportation}. We end this section by extending the former two results in the following theorem.
\begin{theorem}
\label{thm:hmm:t1}
If the $r$-step transition matrix $Q^{r}_{\theta}$ of $(X_{t})_{t\geq 1}$ is positive entrywise for some $r >0$ and $\theta \in \Theta$, then
\begin{align}
\label{mixing}
D_{\theta} = \frac{1}{2} \sum_{t=1}^{\infty}\sup_{a,a'} \sum_{b \in \mathcal{X}}|Q_{\theta}^{t}(a,b)-Q_{\theta}^{t}(a',b)| < \infty,
\end{align}
and $\mathbb{P}_{\theta}^{(n)}(Y_{1}^{n} \in \cdot) \in T_{1}(n C_{H}^{\theta})$ with respect to $(\mathcal{Y}^{n},d_{l_{1}})$, where $C_{H}^{\theta}$ is given as:
\begin{enumerate}
\item In the case $(\mathcal{Y},d_{\mathcal{Y}})$ is countable, $C_{H}^{\theta} = (D_{\theta}+1)^{2};$
\item In the case $(\mathcal{Y},d_{\mathcal{Y}})$ is uncountable, if for any $a \in \mathcal{X}$, the emission distribution $\mathbb{P}_{\theta}(Y_{1}\in \cdot |X_{1} = a) \in T_{1}(C_{\mathcal{Y}})$ for some constant $C_{\mathcal{Y}}$,
and there exists a constant $L>0$ such that for any $a,b \in \mathcal{X}$,
$$W_{1}^{d_{\mathcal{Y}}}\{ \mathbb{P}_{\theta}(Y_{1} \in \cdot \mid X_{1} = a) ,\mathbb{P}_{\theta}(Y_{1} \in \cdot \mid X_{1} = b) \} \leq L 1_{a\neq b},$$
then $C_{H}^{\theta} = C_{\mathcal{Y}} + L^{2}(D_{\theta}+1)^{2}.$
\end{enumerate}
\end{theorem}

\section{Main results}
\label{sec:main}

Consider the setup for studying asymptotic properties of the posterior distribution of $\theta$. Assume that the observe data $Y_{1}^{n}$ has $\PP_{\theta_0}^{(n)}$ as its distribution. Given a prior distribution $\Pi_{n}$ on the parameter space $\Theta$ with respect to Borel $\sigma$-field $\mathcal{F}$, the posterior distribution of $\theta$ conditional on $Y_1^n$ is
\begin{align}
\label{eq:post}
\Pi_{n}(B \mid Y_1^n) = \frac{\int_{B} d\PP^{(n)}_{\theta}(Y_1^n) \, d\Pi_{n}(\theta)}{\int_{\Theta} d\PP^{(n)}_{\theta}(Y_{1}^{n}) \, d\Pi_{n}(\theta)},\quad B \in \mathcal{F}.
\end{align}
The posterior distribution is consistent if it concentrates on arbitrary small neighborhoods of $\theta_0$ with $\PP_{\theta_0}^{(n)}$ probability tending to 1 as $n$ tends to infinity. Let $\| \cdot \|_2$ be the Euclidean distance. Then, the contraction rate in metric space $(\Theta, \|\cdot\|_2)$ is at least $\epsilon_n$ if
\begin{align*}
\Pi_{n}\{\theta \in \Theta: \|\theta - \theta_0\|_2 > M_n \epsilon_n \mid Y_1^n \}\rightarrow 0\text{ in $\PP_{\theta_0}^{(n)}$-probability for every $M_n \rightarrow \infty$}.
\end{align*}
The contraction rate measures the size of small neighborhoods of $\theta_0$ on which the posterior distribution puts almost all mass.

A general technique for deriving posterior contraction rates is based on proving the following testing condition \citep[Theorem 7.3]{ghosal2000convergence}. Let $\epsilon_n $ be a sequence satisfying $\epsilon_{n} \rightarrow 0$ and $(n\epsilon_{n}^2)^{-1} = O(1)$, the testing condition says that for a universal constants $ K >0$ and every sufficiently large $j>0$, there exists a sequence test functions $\phi_{n}$ such that
\begin{align}
\label{eq:opt:test}
\EE_{\theta_0}^{(n)}(\phi_n) \leq e^{-Kj^{2}n\epsilon_{n}^2},\quad \sup_{\theta \in \Theta:j \epsilon_{n} <\|\theta - \theta_0\|_{2}\leq 2j \epsilon_{n}} \EE_\theta^{(n)}(1 - \phi_n) \leq e^{-Kj^{2}n\epsilon_{n}^2}.
\end{align}
Once the testing condition with rate $\epsilon_{n}$ in \eqref{eq:opt:test}  is established, the posterior distribution is consistent with the same rate $\epsilon_n$ for any prior distributions $\Pi_{n}$ putting enough amount of mass near the true parameter $\theta_0$ \citep{ghosal2000convergence,GhoVan07}. Furthermore, if the local asymptotic normality property (LAN) of $\{\PP_{\theta}^{(n)}:\theta \in \Theta\}$ also holds, then the Bernstein-von Mises theorem follows immediately \citep{Van00}.

Our main contribution is to construct $\phi_{n}$ in an HMM and prove that it satisfies the testing condition in \eqref{eq:opt:test}. Starting from the simple hypothesis, we construct the likelihood ratio test with exponentially decaying error based on \cite{hu2011transportation}. Applying the $L^1$-transportation cost information inequality for HMMs, we show that the log likelihood function is a Lipschitz function of the extended HMM sequence $(Y_{t},p^{\theta}_{t})_{t = 1}^{n}$, and it concentrates on its mean with sub-Gaussian tail behavior.
Then, for any $\epsilon >0$, some $0<\xi <1$, and any $\theta_1$ with $\|\theta_1 - \theta_0\|_{2} > \epsilon$, we consider hypothesis of $\PP_{\theta_0}^{(n)}$ versus the complement $\{\PP_{\theta}^{(n)}: \theta \in \Theta,\|\theta - \theta_1\|_{2} \leq \xi \epsilon\}$.
By taking the liklihood ratio tests of $\PP_{\theta_0}^{(n)}$ versus $\PP_{\theta'}^{n}$ where $\theta'$ is the center of $\{ \theta \in \Theta,\|\theta - \theta_1\|_{2} \leq \xi \epsilon\}$, the constructed test has an exponential decaying error. Finally, we prove the testing condition \eqref{eq:opt:test} by an entropy bound of $\{\PP_{\theta}^{(n)}: \theta \in \Theta,j \epsilon_{n} <\|\theta - \theta_0\|_{2}\leq 2j \epsilon_{n}\}$.

We require the following assumptions to ensure that a sequence of test functions $\phi_{n}$ exists that satisfies \eqref{eq:opt:test}.
\begin{enumerate}[label=(\subscript{A}{{\arabic*}})]
\item
\label{a1}
For any $\theta \in \Theta$, $\mathbb{P}_{\theta}^{(n)}(Y_{1}^{n} \in \cdot) \in T_{1}( nC_{H})$ on metric space $(\mathcal{Y}^{n},d_{l^{1}})$ for some constant $C_{H}>0$.
\item
\label{a2}
For any $\theta \in \Theta$, there exists constants $\delta_1,\delta_2 >0$ such that the hidden Markov chain $(X_{t})_{t\geq 1}$ is ergodic, $\|\partial_{y}f^{\theta}_{t}\|_{\infty} < \delta_1
$, and $\sum_{t=0}^{\infty}\|\partial_{p}f^{\theta}_{t}\|_{\infty} < \delta_2$, where
\begin{align}
&\|\partial_{y}f^{\theta}_{t}\|_{\infty} = \sup_{y \neq y',p \in \mathcal{E}}\frac{\| f_{t}^{\theta}(y,p)-f_{t}^{\theta}(y',p)\|_{\text{TV}}}{d_{\fY}(y,y')},
\nonumber
\\
&
\|\partial_{p}f_{t}^{\theta}\|_{\infty} = \sup_{p_1\neq p_1' \in \mathcal{E},y_{1}^{t+1} \in \mathcal{Y}^{t+1}} \frac{\|f^{\theta}_{t}(y_{1}^{t+1},p_{1})-f_{t}^{\theta}(y_{1}^{t+1},p_{1}')\|_{\text{TV}}} {\|p_1-p_1'\|_{\text{TV}}}.
\end{align}

\item
\label{a3}
For any $\theta \in \Theta$, the function $\log \left\{\sum_{x_{1}}p_{1}^{\theta}(x_{1})g_{\theta}(Y_{1}|x_{1})\right\}$ is Lipschitz with norm $L/2$ on metric space $(\mathcal{Z} = \mathcal{Y} \times \mathcal{E},d_{\mathcal{Z}} = d_{\mathcal{Y}} + \|\cdot\|_{\text{TV}})$.
\item
\label{a4}
$\Theta \subset \mathbb{R}^{d}$ is compact and $\theta_0$ lies in its interior.

\item
\label{a6}
As $n \rightarrow \infty$, $\frac{1}{n}\mathbb{E}^{(n)}_{\theta_1}\{l_{n}(\theta_1,Y_{1}^{n}) - l_{n}(\theta_{2},Y_{1}^{n})\} \rightarrow J(\theta_1\mid\theta_2)$ uniformly for $\theta_1,\theta_2 \in \Theta$, where
\begin{align*}
J(\theta_1\mid\theta_2) = \int_{\mathcal{Y}\times \mathcal{E}} \log \frac{\sum_{x \in \mathcal{X}}p(x)g_{\theta_1}(y\mid x) }{\sum_{x \in \mathcal{X}}p(x)g_{\theta_2}(y\mid x)} R_{\theta_{1}}(dy,dp),
\end{align*}
and $ R_{\theta_{1}}$ is the stationary distribution of $(Y_{t},p_{t})$ under $\theta_1$.
\item
\label{a7}
There exist constants $\kappa_1, \kappa_2$ such that $0 < \kappa_1 \leq \kappa_2 \leq 2\kappa_1$ and
for any $\theta_1,\theta_2 \in \Theta$, $$ \kappa_1 \|\theta_1 - \theta_2\|_{2} \leq J(\theta_1|\theta_2) \leq \kappa_2 \|\theta_1 - \theta_2\|_{2}.$$
\end{enumerate}

Assumptions \ref{a1}-\ref{a3} are based on the conditions in \citet[Theorem 3.1]{le2000exponential} and are global over $\Theta$. Assumption \ref{a1} states that the probability measure $\PP_{\theta}^{(n)}$ satisfies the $L^{1}$-transportation cost information inequality on $\Theta$ with $nC_{H}$. Theorem \ref{thm:hmm:t1} shows that the constant $C_{H}$ is expressed in terms of the emission distribution and mixing property of hidden Markov chain. Assumption \ref{a2} implies that the extended HMM sequence $(Y_{t},p_{t})_{t=1}^{n}$ is stationary and satisfies $L^{1}$-transportation cost information inequality on $\Theta$ \citep{hu2011transportation}. Assumption \ref{a3} states that the $\log$ likelihood function is Lipschitz on metric space $(\mathcal{Z} ,d_{\mathcal{Z}})$, which is mild and satisfied in many cases \citep[Example 3.1]{le2000exponential}. Assumption \ref{a6} implies the uniform convergence of normalized log likelihood ratio function, which is slightly stronger than \cite[Proposition 4]{leroux1992maximum}. Assumption \ref{a7} specifies the equivalence of $J(\theta_1 \mid \theta_2)$ and $\|\theta_1 - \theta_2\|_{2}$ which is a regular condition in parametric models \citep[Example 7.1]{ghosal2000convergence}.

\begin{theorem}[Testing Condition]
\label{thm:opt:test}
Let $\epsilon_n >0$ with $\epsilon_{n} \rightarrow 0$ such that $(n\epsilon_{n}^{2})^{-1} = O(1)$. Under Assumptions \ref{a1}-\ref{a7}, there exist positive constants $C', K$ and a sequence of test functions $\phi_{n}$ such that for every sufficiently large $j$,
\begin{align}
\label{eq:thm:test}
\mathbb{E}_{\theta_0}^{(n)}(\phi_{n}) \leq C' e^{- Kn \epsilon_{n}^{2}j^{2} },\quad \sup_{\theta \in \Theta: j \epsilon_{n} <\|\theta - \theta_0\|_{2}\leq 2j \epsilon_{n}}\mathbb{E}_{\theta}^{(n)}(1-\phi_{n}) \leq C'e^{-Kn\epsilon_{n}^{2}j^{2}}.
\end{align}
\end{theorem}

Theorem \ref{thm:opt:test} shows that under certain regularity conditions, there exists a sequence of test functions $\phi_{n}$ satisfying the testing condition in \eqref{eq:opt:test} for HMMs. The existsnce of such tests plays a crucial role in establishing the convergence of posterior distributions \citep{ghosal2000convergence, GhoVan07,ghosal2017fundamentals}. We can directly apply the testing condition in \eqref{eq:thm:test} to establish the contraction rate and asymptotic normality of posterior distribution.

\begin{corollary}[Posterior Convergence]
\label{thm:conver}
Let $\epsilon_n >0$ with $\epsilon_{n} \rightarrow 0$ such that $(n\epsilon_{n}^{2})^{-1} = O(1)$. For $k >1$, define
\begin{align} \label{eq:setB}
K(p_{\theta_0}^{(n)},p_{\theta}^{(n)}) &= \int p_{\theta_0}^{(n)} \log(p_{\theta_0}^{(n)}/p_{\theta}^{(n)}) d\mu^{n}, \nonumber\\
V_{k,0}(p_{\theta_0}^{(n)},p_{\theta}^{(n)}) &= \int p_{\theta_0}^{(n)}| \log(p_{\theta_0}^{(n)}/p_{\theta}^{(n)}) - K(p_{\theta_0}^{(n)},p_{\theta}^{(n)}) |^{k}d\mu^{n}, \nonumber\\
B_{n}(\theta_0,\epsilon;k) &= \{\theta \in \Theta: K(p_{\theta_0}^{(n)},p_{\theta}^{(n)}) \leq n\epsilon^{2}, V_{k,0}(p_{\theta_0}^{(n)},p_{\theta}^{(n)}) \leq n^{k/2}\epsilon^{k} \}.
\end{align}
Assume for some $k>1$, the prior distribution $\Pi_{n}$ of $\theta$ satisfies
\begin{align}
\label{cor:2.5}
\frac{\Pi_{n}(\theta \in \Theta: j\epsilon_{n}< \|\theta - \theta_0\|_{2}\leq 2j \epsilon_{n})}{\Pi_{n}(B_{n}(\theta_0,\epsilon_{n};k))} \leq e^{Kn\epsilon_{n}^{2}j^{2}/2}.
\end{align}
Then under the condition of Theorem \ref{thm:opt:test}, for every $M_{n} \rightarrow \infty$, we have that
\begin{align}
\label{thm:conver:1}
\mathbb{E}_{\theta_0}^{(n)}\Pi_{n}(\theta \in \Theta: \|\theta - \theta_0\|_{2} \geq M_{n}\epsilon_{n} \mid Y_{1}^{n}) \rightarrow 0. 
\end{align}
\end{corollary}

\begin{corollary}[Bernstein-von Mises theorem]
\label{thm:bvm}
Given observation $Y_{1}^{n}$, denote the log likelihood as $l_{n}(\theta)$ and $l_{n}'(\theta)$ as the derivative of $l_{n}$ with respect to $\theta$. Assume the HMM model $\{\mathbb{P}_{\theta}^{(n)}:\theta \in \Theta\}$ satisfies the local asymptotic normality condition (LAN):
\begin{align}
\label{eq:lan}
l_{n}\left(\theta_0 + \frac{h}{\sqrt{n}} \right) - l_{n}(\theta_0) = \frac{1}{\sqrt{n}} h^{\top} l'_{n}(\theta_0) - \frac{1}{2} h^{\top}I_{\theta_0}h + o_{\PP_{\theta_0}^{(n)}}(1),
\end{align}
where $h \in \{h\in \mathbb{R}^{d}: \theta_{0} + \frac{h}{\sqrt{n}} \in \Theta\}$ and $l_{n}'(\theta_0)$ weakly converges to $N_{d}(0,I_{\theta_0})$ under $\PP_{\theta_0}^{(n)}$ probability.
Assume the prior distribution has continuous density $\pi(\theta) > 0 $ on $\Theta$ and the Fisher information matrix $I_{\theta_0}$ at $\theta_0$ is nonsingular. Then under the condition of Theorem \ref{thm:opt:test},
\begin{align}
\label{thm:bvm:1}
\mathbb{E}_{\theta_0}^{(n)}\|\Pi_{n}\left\{ \sqrt{n}(\theta - \theta_0) \mid Y_{1}^{n}\right\} - N_{d}(\Delta_{n,0},I_{\theta_0}^{-1})\|_{TV} \rightarrow 0,
\end{align}
where $\Delta_{n,0} = \frac{1}{\sqrt{n}} I_{\theta_0}^{-1} l'_{n}(\theta_0).$
\end{corollary}

\citet{bickel1998asymptotic} have developed the first set of results showing the consistency and asymptotic normality of maximum likelihood estimators of the parameters in an HMM. 
Corollary \ref{thm:bvm} provides similar theoretical gaurantees about Bayesian estimation in HMMs. Specifically, it states that the posterior distribution of $\theta$ given the observations $Y_{1}^{n}$ in an HMM  is asymptotically normal with random mean vector $\theta_0 + \frac{1}{\sqrt{n}}\Delta_{n,0}$ and covariance matrix $(nI_{\theta_0})^{-1}$. Compared to the similar results in \cite{de2008asymptotic}, which are based on Taylor expansion of the log likelihood function, Corollary \ref{thm:bvm} provides an alternative approach to prove the Bernstein-von Mises theorem using the LAN condition \eqref{eq:lan} and testing condition \eqref{eq:thm:test}. The LAN condition is a general regularity assumption on the class of parametric models $\{\PP_{\theta}^{(n)}:\theta \in \Theta\}$, which is satisfied under the theoretical setup in \cite{bickel1998asymptotic,de2008asymptotic}. In our proof, we construct the sequence of test functions based on a likelihood ratio test and prove the testing condition using the optimal transportation inequality for HMMs. 

\section{Proofs}
\label{sec:proof}

In this section, for simplicity, we denote $\mathbb{P}_{\theta}^{(n)}$ and $\mathbb{E}_{\theta}^{(n)}$ as $\mathbb{P}_{\theta}$ and $\mathbb{E}_{\theta}$, respectively.

\subsection*{Proof of Theorem \ref{thm:hmm:t1}}

In the case $\mathcal{Y}$ is countable, the proof is based on \citet[Theorem 1.3]{kontorovich2008concentration}.
\begin{theorem} \label{thm:kc1.1}[\citet[Theorem 1.3]{kontorovich2008concentration}] Suppose $\mathcal{X}$ is a countable space, $\mathcal{F}$ is the collection of all subsets of $\mathcal{X}^{n}$, $\mathbb{P}$ is the probability measure on $(\mathcal{X}^{n},\mathcal{F})$, and $\mathbb{E}$ is the expectation under $\mathbb{P}$. Suppose $\phi: \mathcal{X}^{n} \mapsto \mathbb{R}$ is a $c$-Lipschitz function (with respect to Hamming metric) for some $c>0$, then for any $r >0$,
$$\mathbb{P}\{|\phi - \mathbb{E}\phi| > r\} \leq 2 \exp(-\frac{r^{2}}{2nc^{2}\|\Delta_{n}\|_{\infty}^{2}}).$$
\end{theorem}
Under $\PP_\theta$, $(X_{t})_{t\geq 1}$ is a Markov chain with transition matrix $Q_{\theta}$, and the mixing coefficient $\|\Delta_{n}\|_{\infty}$ is bounded by the Markov contraction coefficients \citep[Lemma 7.1]{kontorovich2008concentration}:
\begin{align}
\label{Martingale}
\|\Delta_n\|_{\infty} \leq \frac{1}{2} \sum_{t=1}^{\infty}\sup_{a,a'} \sum_{b \in \mathcal{X}}|Q_{\theta}^{t}(a,b)-Q_{\theta}^{t}(a',b)|+1 = D_{\theta}+1.
\end{align}
The transportation inequality of HMM with countable $\mathcal{Y}$ is given as an example in \cite[Section 7.2]{kontorovich2008concentration} and $C_{H}^{\theta} = (D_{\theta}+1)^{2}$.

In the case $\mathcal{Y}$ is uncountable, Theorem \ref{thm:ti} implies that proving the second part of Theorem \ref{thm:hmm:t1} is equivalent to showing that for any Lipschitz function $F: \mathcal{Y}^{n} \mapsto \mathbb{R}$ with $\|F\|_{\text{Lip}} \leq \alpha$ (with respect to $d_{l^1}$),
$$\mathbb{E_\theta} e^{\lambda (F - \mathbb{E_\theta}F)} \leq \exp \left\{ \frac{(\lambda\alpha)^{2}}{2} n(C_{\mathcal{Y}} + L^{2}(D_{\theta} + 1)^{2}) \right\},\quad \forall \lambda \in \mathbb{R}.$$
Let $G_{n}(X_{1}^{n}) = \mathbb{E_{\theta}}[F \mid X_{1}^{n}]$ for $n \geq 1$. Then, for $x_{1}^{n} \neq \tilde{x}_{1}^{n}$,
\begin{align*}
|G_{n}(x_{1}^{n}) - G_{n}(\tilde{x}_{1}^{n})| & = |\int F(y_{1}^{n}) d\mathbb{P}_\theta(Y_{1}^{n} = y_{1}^{n}\mid X_{1}^{n} = x_{1}^{n}) - \int F(y_{1}^{n}) d\mathbb{P}_{\theta}(Y_{1}^{n} = y_{1}^{n}\mid X_{1}^{n} = \tilde x_{1}^{n}) |
\\&
=|\int F(y_{1}^{n}) \left\{ d\mathbb{P}_\theta(Y_{1}^{n} = y_{1}^{n}\mid X_{1}^{n} = x_{1}^{n})- d\mathbb{P}_{\theta}(Y_{1}^{n} = y_{1}^{n}\mid X_{1}^{n} = \tilde x_{1}^{n})\right\} |
\\&
\overset{(i)}{\leq}
\sup_{\|\Phi\|_{\text{Lip}}\leq \alpha} \int \Phi(y_{1}^{n}) \left\{ d\mathbb{P}_\theta(Y_{1}^{n} = y_{1}^{n}\mid X_{1}^{n} = x_{1}^{n})- d\mathbb{P}_{\theta}(Y_{1}^{n} = y_{1}^{n}\mid X_{1}^{n} = \tilde x_{1}^{n})\right\}
\\
& \overset{(ii)}{=} \alpha W_{1}^{d_{l_1}} \{ \mathbb{P}_{\theta}(Y_{1}^{n} \in \cdot \mid X_{1}^{n} = x_{1}^{n}) , \mathbb{P}_{\theta}(Y_{1}^{n} \in \cdot \mid X_{1}^{n} = \tilde x_{1}^{n})\} \\
&\overset{(iii)}{=} \alpha \inf _{\Pi} \int \sum_{t=1}^{n} d_{\mathcal{Y}}(y_{t},\tilde{y}_{t}) d \Pi(y_{1}^{n},\tilde y_{1}^{n})
\\
& = \alpha\sum_{t=1}^{n} 1_{x_{t} \neq \tilde{x}_{t}} W_{1}^{d_{\mathcal{Y}}} \{ \mathbb{P}_\theta(Y_{t} \in \cdot | X_{t} = x_t) ,\mathbb{P}_\theta(Y_{k} \in \cdot | X_{k} = \tilde{x}_t) \} \\
&\overset{(iv)}{ \leq} \alpha L \sum_{k=1}^{n}1_{x_{k} \neq \tilde{x}_{k}} = \alpha Ld_{H}(x, \tilde x),
\end{align*}
where $\Phi$ on the RHS of $(i)$ is an $\alpha$-Lipschitz function on $(\mathcal{Y}^{n},d_{l_{1}})$, $(ii)$ follows from the duality form of $W_{1}$ distance, $(iii)$ follows from the definition of $W_{1}^{d_{l_1}}$ where $\Pi$ is the joint distribution with marginal distributions $\mathbb{P}_{\theta}(Y_{1}^{n} \in \cdot |x_{1}^{n})$ and $\mathbb{P}_{\theta}(Y_{1}^{n} \in \cdot |\tilde{x}_{1}^{n})$, and $(iv)$ follows from assumption that $W_{1}^{d_{\mathcal{Y}}}\{ \mathbb{P}_{\theta}(Y_{1} \in \cdot \mid X_{1} = a) ,\mathbb{P}_{\theta}(Y_{1} \in \cdot \mid X_{1} = b) \} \leq L 1_{a\neq b}$;
therefore, $G_{n}$ is a Lipschitz function on $\mathcal{X}^{n}$ with respect to Hamming metric $d_{H}$. Using (\ref{Martingale}), we have that $\mathbb{P}_{\theta}(X_{1}^{n} \in \cdot) \in T_{1}(n(D_{\theta}+1)^{2})$ on $(\mathcal{X}^{n},d_{l_{1}})$. By Theorem \ref{thm:ti}, for any $\lambda >0$,
$$\mathbb{E}_\theta e^{\lambda(G_{n} - \mathbb{E}_\theta G_{n})} \leq \exp \left(\frac{(\lambda\alpha L)^{2}}{2} n(D_{\theta} + 1)^{2} \right).$$
Using independent tensorization, we have $\mathbb{P}_{\theta}(Y_{1}^{n} \in \cdot |X_{1}^{n}) \in T_{1}(nC_{\mathcal{Y}} )$ on $(\mathcal{Y}^{n},d_{H})$, which implies that
$$\mathbb{E}_{\theta} \{ e^{\lambda (F - E_{\theta}[F \mid X_{1}^{n}])}\mid X_{1}^{n}\} = \mathbb{E}_{\theta} \{ e^{\lambda (F - G_{n})}\mid X_{1}^{n}\}\leq \exp\left(\frac{(\lambda\alpha)^{2}}{2}nC_{\mathcal{Y}} \right). $$
Therefore,
\begin{align*}
\mathbb{E}_{\theta} e^{\lambda (F - \mathbb{E}_{\theta}F}) & = \mathbb{E}_{\theta}\left( \mathbb{E}_{\theta} \{ e^{\lambda (F - G_{n})}\mid X_{1}^{n}\} e^{\lambda (G_{n} - \mathbb{E}_{\theta}F)} \right ) \\
& \leq \mathbb{E}_{\theta}\left( \mathbb{E}_{\theta} \left\{ e^{\lambda (F - E_{\theta}[F \mid X_{1}^{n}])}\mid X_{1}^{n} \right \}\right)\exp \left\{\frac{(\lambda\alpha L)^{2}}{2} n (D_{\theta}+1)^{2}\right\} \\
& \leq \exp \left \{\frac{(\lambda\alpha)^{2}}{2} n(C_{\mathcal{Y}} + L^{2}(D_{\theta}+1)^{2}) \right \} .
\end{align*}
\hfill $\square$

\section*{Proof of Theorem \ref{thm:opt:test}}

\noindent We start by proving the following lemma that is based on \citet[Theorem 3.3]{hu2011transportation}.
\begin{lemma}
\label{lemma:simply}
If Assumptions \ref{a1}-\ref{a3} hold, then for $\theta_1 \neq \theta_2 \in \Theta$, $$F_{n}(\theta_1,\theta_2) = l_{n}(\theta_1,Y_{1}^{n}) - l_{n}(\theta_2,Y_{1}^{n}) $$ is a Lipschitz function on $(\mathcal{Z}^{n},d_{l^{1}})$ with $\|F_{n}(\theta_1,\theta_2)\|_{\text{Lip}} \leq L$ and for any $\theta \in \Theta$,
$$\mathbb{E}_{\theta} e^{\lambda \left\{F_{n}(\theta_1,\theta_2) - \mathbb{E}_{\theta}F_{n}(\theta_1,\theta_2) \right\}}
\leq
\exp\left\{\frac{\lambda^{2}}{2}L^{2}nC_{H}(1 + \delta)^{4}\right\},$$
where $\delta = \max(\delta_1,\delta_2)$ in Assumption \ref{a2}.
\end{lemma}
\subsection*{Proof of Lemma \ref{lemma:simply}}
First we show that for any $\theta \in \Theta$, there exists a constant $C_{E}>0$ such that the probability measure induced by the extended HMM $\{Y_{t},p_{t}^{\theta}\}_{t=1}^{n}$ satisfies $T_{1}(nC_{E})$-inequality on metric space $(\mathcal{Z}^{n} = (\mathcal{Y}\times \mathcal{E})^{n}, d_{l_{1}})$, where the $l_1$-metric on $\mathcal{Z}^{n}$ is defined as:
$d_{l_{1}}(z_{1}^{n},\tilde{z}_{1}^{n}) = \sum _{k=1}^{n}d_{\mathcal{Z}}(z_{k},\tilde{z}_{k})$ and
$$d_{\mathcal{Z}}(z_{k},\tilde{z}_{k}) = d_{\mathcal{Y}}(y_{k}, \tilde y_{k}) + \|p_{k}^{\theta} - p_{k}^{\theta}\|_{\text{TV}}, \quad z_{k} = (y_{k},p_{k}^{\theta}),\tilde z_{k} = (\tilde y_{k},\tilde p_{k}^{\theta}) \in \mathcal{Z}.$$
That is, given $p_{1}^{\theta} = r_{\theta} \in \mathcal{E}$ and for any $\theta \in \Theta$, we want to show the probability measure $\mathbb{P}_{\theta}((Y_{i},p_{i}^{\theta})_{i=1}^{n} \in \cdot)\in T_{1}(nC_{E})$ on $(\mathcal{Z}^{n},d_{l^{1}})$, and $C_{E}$ is given as
$$C_{E} = C_{H}(1 + \delta )^{4},\quad \delta = \max(\delta_1,\delta_2).$$

Let $F$ be a Lipschitz function on $(\mathcal{Z}^{n},d_{l_{1}})$ with $\|F\|_{\text{Lip}} \leq 1$. By the Lipschitz property and definition of $d_{l_{1}}$, for $1\leq k \leq n$ we have that,
\begin{align}
\label{eq:lipF}
& \|\partial_{y_{k}}F\|_{\infty} := \sup_{z_{1}^{k-1},y_{k}\neq \tilde{y}_{k},p_{k},z_{k+1}^{n}} \frac{|F(z_{1}^{k-1},(y_{k},p_{k}),z_{k+1}^{n}) - F(z_{1}^{k-1},(\tilde{y}_{k},p_{k}),z_{k+1}^{n})|}{d_{\mathcal{Y}}(y_{k},\tilde{y}_{k})} \leq 1,
\nonumber
\\
& \|\partial_{p_{k}}F\|_{\infty} := \sup_{z_{1}^{k-1},p_{k}\neq \tilde{p}_{k},y_{k},z_{k+1}^{n}} \frac{|F(z_{1}^{k-1},(y_{k},p_{k}),z_{k+1}^{n}) - F(z_{1}^{k-1},(y_{k},\tilde{p}_{k}),z_{k+1}^{n})|}{\|p_{k} - \tilde{p}_{k}\|_{\text{TV}}} \leq 1.
\end{align}
For $k = 2,\ldots,n$, \eqref{eq:funcf} implies that the prediction filters $p_{k}^{\theta}$ only depend on $y_{1}^{k-1}$, so $F(z_{1}^{n})$ can be written as a function $G$ on $\mathcal{Y}^{n}$.

Next, we get a bound for $\|\partial_{y_{k}}G\|_{\infty}$. Suppose $y_{1}^{n} \neq \tilde{y}_{1}^{n}$ only differ at the $k$th coordinates $y_{k} \neq \tilde{y}_{k}$. Then,
$$G(y_{1}^{n}) - G(\tilde{y}_{1}^{n}) =
F(y_{1}^{n},p_{1}^{n}) - F(\tilde{y}_{1}^{n},\tilde{p}_{1}^{n}).$$
For $j \leq k$, $p_{j} = \tilde{p}_{j}$ because of $y_{1}^{k-1} = \tilde y_{1}^{k-1}$. When $j = k+1$, Assumption \ref{a2} implies that
\begin{align}
\label{eq:k+1}
\| p_{k+1} - \tilde{p}_{k+1}\|_{\text{TV}} &= \| f^{\theta}(y_{k},p_{1}^{k}) - f^{\theta}(\tilde{y_{k}},p_{1}^{k})\|_{\text{TV}} \nonumber\\
&\leq \delta_1 d_{\mathcal{Y}}(y_{k},\tilde{y}_{k});
\end{align}
Finally, for $j \geq k+2$,
\begin{align}
\label{eq:k+2}
\|p_{j} - \tilde{p}_{j}\|_{\text{TV}} &= \|f^{\theta}_{j-k-2}(y_{k+1}^{j-1},p_{k+1}) - f^{\theta}_{j-k-2}(y_{k+1}^{j-1},\tilde{p}_{k+1})\|_{\text{TV}} \nonumber\\
&\leq \|\partial_{p_0}f^{\theta}_{j-k-2}\|_{\infty} \cdot \|p_{k+1} - \tilde{p}_{k+1}\|_{\text{TV}};
\end{align}
therefore, if $y_{1}^{n} , \tilde{y}_{1}^{n}$ differ only at the $k$th coordinates, then
\begin{align}
|G(y_{1}^{n}) - G(\tilde{y}_{1}^{n})| &= | F(y_{1}^{n},p_{1}^{n}) - F(\tilde{y}_{1}^{n},\tilde{p}_{1}^{n})| \nonumber\\
&\leq \| \partial_{y_k} F \|_{\infty}d_{\mathcal{Y}}(y_{k},\tilde{y}_{k}) + \sum_{j=k+1}^{n} \|\partial_{p_j}F\|_{\infty}\|p_{j} - \tilde{p}_{j}\|_{\text{TV}}\nonumber\\
&\overset{(i)}{\leq} d_{\mathcal{Y}}(y_{k},\tilde{y}_{k}) + \left(\|\partial_{p_{k+1}}F\|_{\infty} + \sum_{j=k+2}^{n}\|\partial_{p_0}f^{\theta}_{j-k-2}\|_{\infty}\right)\|p_{k+1} - \tilde{p}_{k+1}\|_{\text{TV}}\nonumber \\
&\overset{(ii)}{\leq} d_{\mathcal{Y}}(y_{k},\tilde{y}_{k})\left(1 + \delta_1 + \delta_1 \delta_2\right) < d_{\mathcal{Y}}(y_{k},\tilde{y}_{k})\left(1 + \delta \right)^{2},
\end{align}
where $(i)$ follows from \eqref{eq:lipF} and $(ii)$ follows from \eqref{eq:k+1}, \eqref{eq:k+2}, and Assumption \ref{a2}. Under $\PP_\theta$, $G(y_{1}^{n})$ is a Lipschitz function on $(\mathcal{Y}^{n},d_{l_1})$ with $\|G\|_{\text{Lip}} \leq (1 + \delta)^{2}$.
Using Assumption \ref{a1}, $\mathbb{P}_{\theta}(Y_{1}^{n} \in \cdot) \in T_{1}( nC_{H})$; Theorem \ref{thm:ti} implies that for any $\lambda \in \mathbb{R}$,
\begin{align}
\label{eq:lG}
\mathbb{E}_{\theta} e^{\lambda \left\{F(Y_{1}^{n},p_{1}^{n}) - \mathbb{E}_{\theta}F(Y_{1}^{n},p_{1}^{n})\right\}} &= \mathbb{E}_{\theta} e^{\lambda\{G(Y_{1}^{n}) - \mathbb{E}_{\theta}G(Y_{1}^{n})\} }\nonumber \\
& \leq \exp\left\{\frac{\lambda^{2}}{2}nC_{H}\left(1 + \delta \right)^{4}\right\}.
\end{align}
Hence, $\mathbb{P}_{\theta}((Y_{i},p_{i}^{\theta})_{i=1}^{n} \in \cdot)\in T_{1}(nC_{E})$ with $C_{E} = C_{H}(1 + \delta )^{4}.$

Concluding the proof of the lemma, we show that $F_{n}(\theta_1,\theta_2)$ is a Lipscthiz function on $(\mathcal{E}^{n},d_{l_1})$. Expressing the log-likelihood additively, we have
\begin{align}
F_{n}(\theta_1,\theta_2) & = l_{n}(\theta_1,Y_{1}^{n}) - l_{n}(\theta_2,Y_{1}^{n})
\nonumber
\\&
= \sum_{t=1}^{n}\left( \log \sum_{x_{t}}p_{t}^{\theta_1}(x_{t})g_{\theta_1}(Y_{t}|x_{t}) -
\log \sum_{x_{t}}p_{t}^{\theta_2}(x_{t})g_{\theta_2}(Y_{t}|x_{t}) \right).
\end{align}
By Assumption \ref{a3}, for $i = 1,2$, the function $\log \sum_{x_{t}}p_{t}^{\theta_i}(x_{t})g_{\theta_i}(Y_{t}|x_{t})$ is a Lipschitz with norm $L/2$ on metric space $(\mathcal{Z} = \mathcal{Y} \times \mathcal{E},d_{\mathcal{Z}} = d_{\mathcal{Y}} + \|\cdot\|_{\text{TV}})$. So $F_{n}(\theta_1,\theta_2)$ is a Lipshitz function on $(\mathcal{Z}^{n},d_{l_1})$ with $\|F_{n}\|_{\text{Lip}} \leq L$. Using Theorem \ref{thm:ti}, we have that
$$\mathbb{E}^{n}_{\theta} e^{\lambda \left\{F_{n}(\theta_1,\theta_2) - \mathbb{E}_{\theta}^{n}F_{n}(\theta_1,\theta_2) \right\}}
\leq
\exp\left\{\frac{\lambda^{2}}{2}L^{2}nC_{H}(1 + \delta )^{4}\right\}.$$

\hfill $\square$

As a special case of Lemma \ref{lemma:simply}, we have
\begin{align}
\label{lem:spec}
\mathbb{E}^{n}_{\theta_1} e^{\lambda \left\{F_{n}(\theta_1,\theta_2) - nH_{n}(\theta_1|\theta_2) \right\}}
\leq
\exp\left\{\frac{\lambda^{2}}{2}L^{2}nC_{H}(1 + \delta )^{4}\right\}.
\end{align}
Now we prove Theorem \ref{thm:opt:test} in three steps. First, we show the existence of a test for $\PP_{\theta_0}$ versus $\PP_{\theta_1}$ with exponentially decaying error probabilities. Second, we show the existence of a test for $\PP_{\theta_0}$ versus the complement $\{\PP_{\theta}:\theta \in \Theta, \|\theta - \theta_1\|_{2} \leq \xi \epsilon\}$ with exponentially decaying error probabilities for any $\epsilon >0$, some $\xi <1$, and any $\theta_1 \in \Theta$ satisfying $\|\theta_1 - \theta_0\|_{2} > \epsilon$. Finally, the proof is completed by showing the existence of a test for $\PP_{\theta_0}$ versus the complement $\{\PP_{\theta}:\theta \in \Theta,\|\theta - \theta_0\|_{2} \geq M \epsilon\}$ with exponentially decaying error probabilities for any $\epsilon >0$ and sufficiently large $M>0$.

\vspace{3mm}

\noindent
\emph{Step 1: Show the existence of a test with exponentially decaying error probabilities for simple hypotheses.}

For $\theta_1 \neq \theta_0 \in \Theta$, we consider the simple hypotheses
\begin{align}
\label{hy:sim}
H_{0}: \theta = \theta_0 \quad \text{vs} \quad H_{1}: \theta = \theta_1.
\end{align}
Based on Lemma \ref{lemma:simply}, we show the existence of a test with Type \Romannum{1} and \Romannum{2} error probability decaying exponentially to $0$ as $n\rightarrow \infty$. For $\theta \neq \theta' \in \Theta$, define $$H_{n}(\theta|\theta') = \frac{1}{n}\mathbb{E}_{\theta}\{l_{n}(\theta,Y_{1}^{n}) - l_{n}(\theta',Y_{1}^{n})\}.$$
Using Jensen's inequality, we have $H_{n}(\theta|\theta') \geq 0$. By Assumption \ref{a6}, $H_{n}(\theta | \theta') \rightarrow J(\theta | \theta') >0$ as $n \rightarrow \infty$; so there exists an integer $N = N(\theta,\theta')$ depending on $(\theta, \theta')$ such that $H_{n}(\theta | \theta') \geq \frac{1}{2}J(\theta| \theta')$ for any $n > N(\theta,\theta')$. Let the critical value $r \in (0, \frac{1}{4}J(\theta_0|\theta_1)).$ We construct the likelihood ratio test for the hypotheses in \eqref{hy:sim} as
\begin{align}
\label{hy:sim:test}
\phi_{\theta_1} = 1\{l_{n}(\theta_0,Y_{1}^{n}) - l_{n}(\theta_1,Y_{1}^{n}) \leq nr\} \equiv 1\{F_{n}(\theta_0,\theta_1) \leq nr\},
\end{align}
and reject $H_0$ if $\phi_{\theta_1} = 1$.
For any $n > N(\theta_0,\theta_1)$, we have $ H_{n}(\theta_0 | \theta_1) - r >\frac{1}{4}J(\theta_0 | \theta_1)$ and the probability of Type \Romannum{1} error of the test in \eqref{hy:sim:test} is
\begin{align}
\label{eq:sim:e1}
\mathbb{P}_{\theta_0}(\text{reject} \;H_0) &= \mathbb{P}_{\theta_0}(\phi_{\theta_1}=1)= \PP_{\theta_0}[l_{n}(\theta_0,Y_{1}^{n}) - l_{n}(\theta_1,Y_{1}^{n}) \leq n r ]
\nonumber
\\&
= \PP_{\theta_0}[l_{n}(\theta_0,Y_{1}^{n}) - l_{n}(\theta_1,Y_{1}^{n}) - nH_{n}(\theta_0|\theta_1) \leq n\{r - H_{n}(\theta_0|\theta_1)\}]
\nonumber
\\&
= \PP_{\theta_0}[F_{n}(\theta_0,\theta_1) - \mathbb{E}_{\theta_0}F_{n}(\theta_0,\theta_1) \leq n\{r - H_{n}(\theta_0|\theta_1)\}]
\nonumber
\\&
\overset{(i)}{\leq} \exp\left(-\frac{n(r - H_{n}(\theta_0|\theta_1))^{2}}{\tilde{C}}\right)
\overset{(ii)}\leq \exp\left(-\frac{nJ(\theta_0 |\theta_1)^{2}}{16 \tilde{C}}\right),
\end{align}
where $(i)$ follows from \eqref{lem:spec} with $\tilde{C} = L^{2}C_{H}(1 + \delta )^{4}$, and $(ii)$ follows from $ H_{n}(\theta_0 | \theta_1) - r >\frac{1}{4}J(\theta_0 | \theta_1) > 0$.

The exponentially decaying bound on the probability of Type \Romannum{2} error follows similarly by switching the roles of $\theta_0$ and $\theta_1$. There exists $N(\theta_1,\theta_0) \in \mathbb{N}$ such that for any $n > N(\theta_1,\theta_0)$, $H_{n}(\theta_1|\theta_0) \geq \frac{1}{2}J(\theta_1|\theta_0)$, and the probability of Type \Romannum{2} error of $\phi_{\theta_1}$ in \eqref{hy:sim:test} is
\begin{align}
\label{eq:sim:e2}
\mathbb{P}_{\theta_1}(\text{fail to reject}\; H_0) &= \mathbb{P}_{\theta_1}(\phi_{\theta_1}=0)=\mathbb{P}_{\theta_1}[l_{n}(\theta_0,Y_{1}^{n}) - l_{n}(\theta_1,Y_{1}^{n}) > nr ]
\nonumber
\\&
= \mathbb{P}_{\theta_1}[l_{n}(\theta_0,Y_{1}^{n}) - l_{n}(\theta_1,Y_{1}^{n}) + nH_{n}(\theta_1|\theta_0) > n\{H_{n}(\theta_1|\theta_0) + r \} ]
\nonumber
\\&
= \mathbb{P}_{\theta_1}[F_{n}(\theta_1,\theta_0) - \mathbb{E}_{\theta_1} F_{n}(\theta_1,\theta_0)< -n\{ H_{n}(\theta_1|\theta_0) + r \} ]
\nonumber
\\&
\leq \exp\left(-\frac{n(H_{n}(\theta_1|\theta_0) + r )^{2}}{\tilde{C}}\right)
\leq \exp\left(-\frac{nJ(\theta_1|\theta_0)^{2}}{4\tilde{C}}\right) .
\end{align}
By choosing $n > \max\{ N(\theta_1,\theta_0), N(\theta_1,\theta_0)\},$ we can see that the probabilities of Type \Romannum{1} and Type \Romannum{2} error of test \eqref{hy:sim:test} decay exponentially in $n$.
\vspace{3mm}

\noindent
\emph{Step 2: For any $\epsilon >0$, some $\xi <1$, any $\theta_1 \in \Theta$ satisfying $\|\theta_1 - \theta_0\|_{2} > \epsilon$, show the existence of a test for $\PP_{\theta_0}$ versus the complement $\{\PP_{\theta}: \theta \in \Theta, \|\theta - \theta_1\|_{2} \leq \xi \epsilon\}$ with exponentially decaying error probabilities.}

For any $\epsilon >0$, let $\theta' \in \Theta$ satisfy $\epsilon<\|\theta' - \theta_0\|_{2} < 2\epsilon$, $U \subset \{\theta \in \Theta, \epsilon<\|\theta' - \theta_0\|_{2} < 2\epsilon\}$ be an open neighborhood with diameter $\frac{\epsilon}{2}$, and $\theta_U \in U$ be a center such that $\|\theta - \theta_U\|_{2}<\frac{\epsilon}{4}$ for any $\theta \in U$. We consider the following hypotheses:
\begin{align}
\label{hy:ball}
H_{0}: \theta = \theta_0 \quad \text{vs} \quad H_{1}: \theta \in U.
\end{align}
We want to show the existence of a test $\phi_{U}^{\epsilon}$ with exponentially decaying Type \Romannum{1} and \Romannum{2} error probabilities; that is,
$$ \mathbb{E}_{\theta_0}(\phi^{\epsilon}_{U}) \lesssim e^{-n\tilde{c}\epsilon^{2}},\quad\sup_{\theta \in U} \, \mathbb{E}_{\theta}(1-\phi_{U}^{\epsilon}) \lesssim e^{-n\tilde{c}\epsilon^{2}},$$
where $\lesssim$ denotes inequality up to a fixed constant. Let $r >0$ be the critical value to be chosen later. Then, define
\begin{align}
\label{hy:ball:test}\phi_{U}^{\epsilon} =1\{l_{n}(\theta_0,Y_{1}^{n}) - l_{n}(\theta_U,Y_{1}^{n}) \leq nr\},
\end{align}
and reject $H_0$ in \eqref{hy:ball} if $\phi_{U}^{\epsilon} =1$.

The corresponding probability of Type \Romannum{1} error of $\phi_{U}^{\epsilon} $ is
\begin{align}
\label{er:b:1}
\mathbb{P}_{\theta_0}(\text{reject $H_0$}) &= \mathbb{E}_{\theta_0}(\phi_{U}^{\epsilon})
\nonumber
\\&
= \mathbb{P}_{\theta_0}[l_{n}(\theta_0,Y_{1}^{n}) - l_{n}(\theta_U,Y_{1}^{n}) \leq nr]
\nonumber
\\&
\overset{(i)}{\leq}\exp\left(-\frac{n(r - H_{n}(\theta_0|\theta_U))^{2}}{\tilde{C}}\right).
\end{align}
where $(i)$ follows from Step 1 provided $r - H_{n}(\theta_0|\theta_U) <0$.

The corresponding probability of Type \Romannum{2} error of $\phi_{U}^{\epsilon} $ is
\begin{align}
\label{er:b:2}
\sup_{\theta \in U} \, & \mathbb{P}_{\theta}(\text{fail to reject $H_0$})
=\sup_{\theta \in U} \mathbb{E}_{\theta}(1 - \phi_{U}^{\epsilon})
\nonumber
\\
&=
\sup_{\theta \in U} \mathbb{P}_{\theta}[l_{n}(\theta_0,Y_{1}^{n}) - l_{n}(\theta_U,Y_{1}^{n}) > nr ]
\nonumber
\\&
=
\sup_{\theta \in U} \mathbb{P}_{\theta}[l_{n}(\theta_0,Y_{1}^{n}) -l_{n}(\theta,Y_{1}^{n}) + l_{n}(\theta,Y_{1}^{n}) - l_{n}(\theta_U,Y_{1}^{n}) > nr]
\nonumber
\\&
=
\sup_{\theta \in U} \mathbb{P}_{\theta}[l_{n}(\theta_0,Y_{1}^{n}) -l_{n}(\theta,Y_{1}^{n}) + nH_{n}(\theta | \theta_0)
\nonumber
\\&
\quad
+ l_{n}(\theta,Y_{1}^{n}) - l_{n}(\theta_U,Y_{1}^{n}) - nH_{n}(\theta|\theta_U) > n(r + H_{n}(\theta | \theta_0) - H_{n}(\theta|\theta_U)) ]
\nonumber
\\&
\leq \sup_{\theta \in U}\{ \mathbb{P}_{\theta}[l_{n}(\theta,Y_{1}^{n}) -l_{n}(\theta_0,Y_{1}^{n}) - nH_{n}(\theta | \theta_0) < -\frac{n}{2}(r + H_{n}(\theta | \theta_0) - H_{n}(\theta|\theta_U)) ]
\nonumber
\\&
\quad
+ \mathbb{P}_{\theta} [l_{n}(\theta,Y_{1}^{n}) - l_{n}(\theta_U,Y_{1}^{n}) - nH_{n}(\theta|\theta_U) > \frac{n}{2}(r + H_{n}(\theta | \theta_0) - H_{n}(\theta|\theta_U)) ]\}
\nonumber
\\&
\overset{(i)}{\leq} 2\sup_{\theta \in U}\left\{ \exp\left(-\frac{n(r + H_{n}(\theta | \theta_0) - H_{n}(\theta|\theta_U) )^{2}}{4\tilde{C}}\right) \right\}.
\end{align}
where $(i)$ follows from Lemma \ref{lemma:simply} provided $r + H_{n}(\theta | \theta_0) - H_{n}(\theta|\theta_U) >0$.

Now we choose the critical value $r$ such that
$$\frac{1}{4}J(\theta_0|\theta_U) <r < \frac{1}{2}J(\theta_0|\theta_U).$$
From Step 1, there exists $N(\theta_0,\theta_U) \in \mathbb{N}$ such that for any $n > N(\theta_0,\theta_U)$, $H_{n}(\theta_0|\theta_U) \geq \frac{3}{4}J(\theta_0|\theta_U)$; so we have that
\begin{align}
\label{er:b:21}
r - H_{n}(\theta_0|\theta_U) \leq -\frac{1}{4} J(\theta_0|\theta_U)
\overset{(i)}{\leq} -\frac{\kappa_2}{4}\|\theta_0 - \theta_U\|_{2}
\leq -\frac{9\kappa_2}{16}\epsilon,
\end{align}
where $(i)$ follows from Assumption \ref{a7}. Combining the results in \eqref{er:b:1} and \eqref{er:b:21} gives the exponentially decaying bounds for the Type \Romannum{1} error:
$$\mathbb{E}_{\theta_0}(\phi_{U}^{\epsilon}) \leq
\exp\left[- n\epsilon^{2} \frac{81\kappa_2^{2}}{256\tilde{C}} \right ] .$$

Assumption \ref{a6} implies that $H_{n}(\theta |\theta_0) \rightarrow J(\theta | \theta_0)$ uniformly for any $\theta \in U$; so there exists $N'(\theta_0) \in \mathbb{N}$ such that for any $n > N'(\theta_0)$, $H_{n}(\theta |\theta_0) \geq \frac{15}{16}J(\theta |\theta_0)$ for any $\theta \in U$. Furthermore, by choosing $n > \max\{N'(\theta_0), N'(\theta_U)\}$, we also have $\frac{15}{16}J(\theta |\theta_U)\leq H_{n}(\theta |\theta_U) \leq \frac{17}{16}J(\theta |\theta_U)$ for any $\theta \in U$. Therefore,
\begin{align}
\label{er:b:22}
\inf_{\theta \in U}\{r &+ H_{n}(\theta | \theta_0) - H_{n}(\theta|\theta_U)\}
\geq \frac{1}{4}J(\theta_0|\theta_U) + \inf_{\theta \in U}\{ \frac{3}{4}J (\theta | \theta_0) - \frac{4}{5}J(\theta |\theta_U)\}
\nonumber
\\
&
\overset{(i)}{\geq} \frac{\kappa_1}{4}\|\theta_0 - \theta_U\|_{2} + \frac{3\kappa_{1} }{4}\inf_{\theta \in U}\|\theta - \theta_0\|_{2} - \frac{4\kappa_{2} }{5}\sup_{\theta \in U} \|\theta - \theta_U\|_{2}
\nonumber
\\
&
\overset{(ii)}{\geq} (\frac{9\kappa_{1} }{16} - \frac{\kappa_2}{5})\epsilon >0,
\end{align}
where $(i)$ follows from Assumption \ref{a7} and $(ii)$ follows from for $\theta \in U$, $\|\theta_0 - \theta_U\|_{2} \geq \frac{3\epsilon}{4}, \|\theta_0 - \theta\|_{2} \geq \frac{\epsilon}{2} $, and $\|\theta - \theta_U\|_{2} \leq \frac{\epsilon}{4}.$
Combining the results in \eqref{er:b:2} and \eqref{er:b:22} gives the exponentially decaying bounds for the Type \Romannum{2} error:
$$\sup_{\theta \in U} \mathbb{E}_{\theta}(1 - \phi_{U}^{\epsilon}) \leq
2\exp\left[- n\epsilon^{2} \frac{(\frac{9\kappa_{1} }{16} - \frac{\kappa_2}{5})^{2}}{4\tilde{C}} \right ] .$$

By choosing $n > \max\{ N(\theta_0,\theta_U),N(\theta_0), N'(\theta_0)\},$ the probabilities of Type \Romannum{1} and Type \Romannum{2} error of test \eqref{hy:ball:test} is dominated by $\exp(-\tilde c n\epsilon^{2})$ where $\tilde c = \min\{\frac{81\kappa_2^{2}}{256\tilde{C}} ,\frac{(\frac{9\kappa_{1} }{16} - \frac{\kappa_2}{5})^{2}}{4\tilde{C}}\}$. Furthermore, this result holds for any open ball $U_{\xi}$ of radius $\xi\epsilon$ such that $ U_{\xi} \subset \{\theta \in \Theta, \epsilon<\|\theta' - \theta_0\|_{2} < 2\epsilon\}$.

\vspace{3mm}
\noindent
\emph{Step 3: For any $\epsilon >0$ and sufficiently large $M>0$, show the existence of a test for testing $\PP_{\theta_0}$ versus the complement $\{\PP_{\theta}: \theta \in \Theta, \|\theta - \theta_0\|_{2} \geq M \epsilon\}$ with exponentially decaying error probabilities.}

We want to show the existence of a test $\phi^{\epsilon}_{M}$ with exponentially decaying Type \Romannum{1} and \Romannum{2} error probabilities; that is,
$$ \mathbb{E}_{\theta_0}(\phi^{\epsilon}_{M}) \lesssim e^{-c_{1}n\epsilon^{2}M^{2}},\quad\sup_{\theta \in \Theta:\|\theta - \theta_0\|_{2} > M/\sqrt{n}}\mathbb{E}_{\theta}(1-\phi^{\epsilon}_{M}) \lesssim e^{-c_{2}n\epsilon^{2}M^{2}}.$$
We start by defining a collection of sets that are used to construct our test later.
For any $\epsilon >0$ and $\xi \in (0,1)$, the compactness of the closure of $\{\theta'\in \Theta, \epsilon < \|\theta' - \theta_0\|_{2} \leq 2\epsilon\}$ implies that there exist a finite number of open balls $\{U_{l,1}\}_{l=1}^{N(\xi, \epsilon)}$ with radius $\xi\epsilon$ such that
$$ \{\theta'\in \Theta, \epsilon < \|\theta' - \theta_0\|_{2} < 2\epsilon\} \subseteq \bigcup_{l=1}^{N(\xi, \epsilon)} U_{l,1},$$
where $N(\xi, \epsilon) = N(\xi\epsilon, \{\theta'\in \Theta, \epsilon < \|\theta' - \theta_0\|_{2} \leq 2\epsilon\},\|\cdot\|_{2})$ is the minimum number of balls with radius $\xi\epsilon$ under metric $\|\cdot\|_{2}$ needed to cover the set $\{\theta'\in \Theta, \epsilon < \|\theta' - \theta_0\|_{2} \leq 2\epsilon\}$. Extending this argument further, we have that for a large number $M$ and any integer $j \geq M$, there are $N(\xi j\epsilon, \{\theta'\in \Theta, j\epsilon < \|\theta' - \theta_0\|_{2} \leq 2j\epsilon\},\|\cdot\|_{2})$ many open balls with radius $\xi j \epsilon$ under metric $\|\cdot\|_{2}$ that cover the set $\{\theta'\in \Theta, j\epsilon < \|\theta' - \theta_0\|_{2} \leq 2j\epsilon\}$.
For ease of notation and following $N(\xi, \epsilon)$ defined earlier, denote $$N(\xi,j\epsilon) = N(\xi j \epsilon, \{\theta'\in \Theta, j\epsilon < \|\theta' - \theta_0\|_{2} \leq 2j\epsilon\},\|\cdot\|_{2})$$
for any $\epsilon >0$, $0< \xi<1$, and $j \geq M$. The compactness of the closure of $ \{\theta'\in \Theta, j\epsilon < \|\theta' - \theta_0\|_{2} \leq 2j\epsilon\}$ for every integer $j \geq M$ implies that we have a collection of open balls $\{U_{l,j}\}_{l=1}^{N(\xi,j\epsilon)}$ with radius $\xi\epsilon$ satisfying
\begin{align*}
\{\theta'\in \Theta, j\epsilon < \|\theta' - \theta_0\|_{2} \leq 2j\epsilon\} \subseteq \bigcup_{l=1}^{N(\xi, j\epsilon)} U_{l,j}, \quad j \geq M.
\end{align*}

We define the test $\phi^{\epsilon}_{M}$ using the test constructed in Step 2 for different choices of the open ball $U$. For each ball in $\{U_{l,j}\}_{l=1}^{N(\xi, j\epsilon)}$, we can find a collection of tests $\{\psi_{U_{l,j}}^{n}\}_{l=1}^{N(\xi, j\epsilon)}$ with Type \Romannum{1} and \Romannum{2} error rates dominated by $\exp(-\tilde{c}nj^{2}\epsilon^{2})$; see Step 2 for the definition of $\tilde c$. Let $\phi_{M}^{\epsilon}$ be the supremum of countably many tests obtained this way:
$$ \phi_{M}^{\epsilon} = \sup_{j\geq M} \sup_{1\leq l \leq N(\xi,j\epsilon)} \psi_{U_{j,l}}^{n}.$$
Then, for a sufficiently large $n$,
\begin{align}
\mathbb{E}_{\theta_0}\phi_{M}^{\epsilon} &\leq \sum_{j=M}^{\infty}\sum_{l=1}^{N(\xi,j\epsilon)} \mathbb{E}_{\theta_0}\psi^{n}_{U_{j,l}}
\overset{(i)}{\leq} \sum_{j=M}^{\infty} N(\xi,j\epsilon) e^{-\tilde{c}nj^{2}\epsilon^{2}}
\leq \sup_{j\geq 1}N(\xi, j \xi) \sum_{j=M}^{\infty} e^{-\tilde{c}nj^{2}\epsilon^{2}}
\nonumber \\&
\leq \frac{\sup_{j\geq 1}N(\xi, j \xi)}{1 - e^{-\tilde c n\epsilon^{2}}} e^{-\tilde{c}nM^{2}\epsilon^{2}},
\label{eq:test:1}
\end{align}
where $(i)$ follows from the upper bound on the Type \Romannum{1} error probability obtained in Step 2.
Similarly, the upper bound on the Type \Romannum{2} error probability obtained in Step 2 implies that, for a sufficiently large $n$,
\begin{align}
\sup_{\theta \in \Theta:\|\theta - \theta_0\|_{2}> M\epsilon}\mathbb{E}_{\theta}(1-\phi_{M}^{\epsilon})
&\overset{(i)}{\leq} \sup_{U_{l,j},\;\;j \geq M} \mathbb{E}_{\theta}(1-\psi_{U_{l,j}}^{n}) \nonumber \\
&\leq \sup_{j \geq M} 2e^{-\tilde{c}nj^{2}\epsilon^{2} } \leq 2 e^{-\tilde{c}nM^{2}\epsilon^{2}},
\label{eq:test:2}
\end{align}
where $(i)$ follows from the construction of $\phi_{M}^{\epsilon}$ because for each $\theta \in \{\theta \in \Theta: \|\theta - \theta_0\|_{2}> M\epsilon\}$, there exists a $j \geq M$ and a test $\psi_{U_{l,j}}^{n}$ with $\phi \geq \psi_{U_{l,j}}^{n}$ satisfying $\mathbb{E}_{\theta}( 1- \psi_{U_{l,j}}^{n}) \leq
e^{-\tilde{c}nj^{2}\epsilon^{2} }.$

Let $\epsilon_{n}>0$, $\epsilon_{n} \rightarrow 0$ and $(n\epsilon_{n}^{2})^{-1} = O(1)$. Then, \citet[Example 7.1]{ghosal2000convergence} implies that
$$\sup_{j\geq 1}N(\xi, j\xi) = \sup_{j\geq 1}N(\xi j \epsilon, \{\theta'\in \Theta, j\epsilon < \|\theta' - \theta_0\|_{2} < 2j\epsilon\},\|\cdot\|_{2}) \leq \left(\frac{12}{\xi}\right)^{d} ;$$
therefore, let $\xi = \frac{1}{4}$; \eqref{eq:test:1} and \eqref{eq:test:2} imply that there exists a constant $$\tilde{c} = \min\{\frac{81\kappa_2^{2}}{256L^{2}C_{H}(1 + \delta )^{4}} ,\frac{(\frac{9\kappa_{1} }{16} - \frac{\kappa_2}{5})^{2}}{4L^{2}C_{H}(1 + \delta )^{4}}\} $$ and a sequence of test functions $\phi_{n} = \phi_{M}^{\epsilon_{n}}$ for large enough $M$ such that for every sufficiently large $n$,
$$ \mathbb{E}_{\theta_0}(\phi_{n}) \leq \frac{\left(48\right)^{d} }{1 - e^{-\tilde{c}}} e^{-\tilde{c}n M^{2}\epsilon_{n}^{2}},\quad\sup_{\theta \in \Theta:\|\theta - \theta_0\| > M\epsilon_{n}}\mathbb{E}_{\theta}(1-\phi_{n}) \leq 2e^{-\tilde{c}nM^{2}\epsilon_{n}^{2}}.$$
The proof is complete by observing that
$$\sup_{\theta \in \Theta: M\epsilon_{n}< \|\theta - \theta_0\|_{2}\leq 2M \epsilon_{n}}\mathbb{E}_{\theta}(1-\phi_{n}) \leq \sup_{\theta \in \Theta:\|\theta - \theta_0\|_{2} > M\epsilon_{n}}\mathbb{E}_{\theta}(1-\phi_{n}) .$$
\hfill $\square$

\section*{Proof of Corollary \ref{thm:conver}}
\noindent
The statement \eqref{thm:conver:1} follows by verifying the conditions in \cite[Theorem 3]{GhoVan07}.
\hfill $\square$

\section*{Proof of Corollary \ref{thm:bvm}}
\noindent The proof is based on \cite[Theorem 10.1]{Van00}. First we verify the contiguity property of $\PP_{\theta}$ and the existence of tests as specified in \citet[Theorem 10.1]{Van00}.

\emph{Contiguity property:} From local asymptotic normality condition \eqref{eq:lan}, under $\PP_0$-probability, we have that for any $h \in \sqrt{n}(\Theta - \theta_0)$
\begin{align}
\label{eq:pf:lan}
l_{n}(\theta_0 + \frac{h}{\sqrt{n}}) - l_{n}(\theta_0) = \frac{1}{\sqrt{n}}h^{\top}l'_{n}(\theta_0) - \frac{1}{2}h^{\top}I_{\theta_0}h + o_{\PP_{\theta_0}}(1),
\end{align}
where $l_{n}'(\theta_0)$ weakly converges to $N_{d}(0,I_{\theta_0})$ under $\PP_{\theta_0}$-probability.
Let $\PP_{\theta_n}$ denote the distribution of $(Y_{1},\ldots,Y_{n})$ under parameter $\theta_n = \theta_0 +\frac{h}{\sqrt{n}}$ for some $h \in \sqrt{n}(\Theta - \theta_0)$.
LeCam's first lemma \citep{lecam1960locally} implies that $\PP_{\theta_n}$ is mutually contiguous to $\PP_{\theta_0}$.

\emph{Existence of test:} If the assumptions of Theorem \ref{thm:opt:test} hold, then for any sequence $\epsilon$ satisfying $\epsilon_{n}>0$, $\epsilon_{n} \rightarrow 0$, and $(n\epsilon_{n}^{2})^{-1} = O(1)$, there exists a constant $K$ and a sequence of test $\phi_{n}$, such that for every sufficiently large $n$,
$$ \mathbb{E}_{\theta_0}(\phi_{n}) \lesssim e^{-KM^{2}n\epsilon_{n}^{2}},\;\quad\sup_{\theta \in \Theta:\|\theta - \theta_0\| > M\epsilon_{n}}\mathbb{E}_{\theta}(1-\phi_{n}) \leq 2e^{-KM^{2}n\epsilon_{n}^{2}}.$$
By setting $\epsilon_{n} = n^{-1/2}$ and $r_{n}$ be any sequence tending to infinity, we can find a sequence of tests $\phi_{n}$ satisfying as $n \rightarrow \infty$,
\begin{align}
\label{eq:bvm:test}
\mathbb{E}_{\theta_0}(\phi_{n}) \lesssim e^{-Kr_{n}^{2}}\rightarrow 0,\quad \sup_{\theta \in \Theta:\|\theta - \theta_0\| > r_{n}/\sqrt{n}}\mathbb{E}_{\theta}(1-\phi_{n}) \leq 2e^{-Kr_{n}^{2}} \rightarrow 0.
\end{align}

Now we prove Corollary \ref{thm:bvm} in two steps. In the first step, we show that the total variance distance between the posterior distribution of $h = \sqrt{n}(\theta - \theta_0)$ relative to prior $\Pi_{n}$ and the restricted prior $\Pi_{n}^{C_{n}}$, where $C_{n}$ is the ball with radius $r_{n}$, is asmptotically negligible as $r_{n} \rightarrow \infty.$
By \eqref{eq:bvm:test}, for every $r_{n} \rightarrow \infty$, we can find a constant $K$ and a sequence of tests $\phi_{n}$, such that as $n \rightarrow \infty$,
\begin{align}
\label{test}
\mathbb{E}_{\theta_0}\phi_{n} \rightarrow 0,\quad\sup_{\theta \in \Theta:\|\theta - \theta_0\| > r_{n}/\sqrt{n}}\mathbb{E}_{\theta}(1 - \phi_{n}) \leq 2e^{-Kr_{n}^{2}}.
\end{align}
Let $\mathcal{F}$ denote the Borel $\sigma$-field on $\sqrt{n}(\Theta - \theta_0)$ and $\Pi_{h|Y_{1}^{n}}^{C_{n}}$ denote the posterior distribution of $h$ conditional on $Y_{1}^{n}$ with respect to the restricted prior $\Pi_{n}^{C_{n}}$. We have that
\begin{align}
\label{eq:bvm:1}
&\|\Pi_{h|Y_{1}^{n}} - \Pi_{h|Y_{1}^{n}}^{C_{n}}\|_{\text{TV}} = \sup_{B \in \mathcal{F}}\left| \Pi_{h|Y_{1}^{n}}(B) - \Pi_{h|Y_{1}^{n}}^{C_{n}}(B)\right|
\nonumber
\\&
= \sup_{B \in \mathcal{F}}\left| \frac{\int 1_{B}(h)\pi(h)p_{\theta + h/\sqrt{n}}(Y_{1}^{n})dh}{\int \pi(s)p_{\theta + s/\sqrt{n}}(Y_{1}^{n})ds} - \frac{\int 1_{C_n}(h) 1_{B}(h)\pi(h)p_{\theta + h/\sqrt{n}}(Y_{1}^{n})dh}{\int 1_{C_n}(s)\pi(s)p_{\theta + s/\sqrt{n}}(Y_{1}^{n})ds} \right|
\nonumber
\\&
= \sup_{B \in \mathcal{F}}\left|\frac{\int \{1_{B}(h) - 1_{B}(h)1_{C_{n}}(h)\}\pi(h)p_{\theta + h/\sqrt{n}}(Y_{1}^{n})dh}{\int \pi(s)p_{\theta + s/\sqrt{n}}(Y_{1}^{n})ds} - \int 1_{C_n}(h) 1_{B}(h)\pi(h)p_{\theta + h/\sqrt{n}}(Y_{1}^{n})dh \times \right.
\nonumber
\\
&
\left.\left[\{\int 1_{C_n}(s)\pi(s)p_{\theta + s/\sqrt{n}}(Y_{1}^{n})ds\}^{-1} - \{\int\pi(s)p_{\theta + s/\sqrt{n}}(Y_{1}^{n})ds\}^{-1}\right]
\right|
\nonumber
\\&
=
\sup_{B \in \mathcal{F}}\left|\frac{\int \{1_{B}(h) - 1_{B}(h)1_{C_{n}}(h)\}\pi(h)p_{\theta + h/\sqrt{n}}(Y_{1}^{n})dh}{\int \pi(s)p_{\theta + s/\sqrt{n}}(Y_{1}^{n})ds} - \right.
\nonumber
\\
&\left.
\frac{\int 1_{C_n}(h) 1_{B}(h)\pi(h)p_{\theta + h/\sqrt{n}}(Y_{1}^{n})dh}{\int 1_{C_n}(s)\pi(s)p_{\theta + s/\sqrt{n}}(Y_{1}^{n})ds }\times \frac{\int\{1 - 1_{C_n}(s)\}\pi(s)p_{\theta + s/\sqrt{n}}(Y_{1}^{n})ds}{\int\pi(s)p_{\theta + s/\sqrt{n}}(Y_{1}^{n})ds} \right|
\nonumber
\\&
= \sup_{B \in \mathcal{F}}| P_{h|Y_{1}^{n}}(B\cap C_{n}^{c}) - P_{h|Y_{1}^{n}}^{C_{n}}(B) P_{h|Y_{1}^{n}}( C_{n}^{c}) |
\nonumber
\\&
\leq 2P_{h|Y_{1}^{n}}(C_{n}^{c})
\end{align}
Let $U$ be a ball with fixed radius. Define $\PP_{n,U}= \int \PP_{\theta + h/\sqrt{n}} d\Pi^{U}_{n}(h)$ and the associated expectation as $\EE_{n,U}$. Using \eqref{eq:lan} and boundedness of $\pi(\cdot)$, $\PP_{n,U}$ is mutually contiguous with $\PP_{\theta_0}$;
so using \eqref{eq:bvm:test}, we have that
\begin{align*}
\mathbb{E}_{n,U} \Pi_{h|Y_{1}^{n}}(C_{n}^{c}) &=\mathbb{E}_{n,U} \Pi_{h|Y_{1}^{n}}(C_{n}^{c})(1 - \phi_{n}) + \mathbb{E}_{0}\Pi_{h|Y_{1}^{n}}(C_{n}^{c})\phi_n + o(1)\\
&= \mathbb{E}_{n,U} \Pi_{h|Y_{1}^{n}}(C_{n}^{c})(1 - \phi_{n}) + o(1),
\end{align*}
where $\mathbb{E}_{n,U} \Pi_{h|Y_{1}^{n}}(C_{n}^{c})(1 - \phi_{n})$ can be expressed as
\begin{align}
\label{eq:bvm:2}
\mathbb{E}_{n,U} \Pi_{h|Y_{1}^{n}}(C_{n}^{c})(1 - \phi_{n})
&= \int \PP_{\theta + h/\sqrt{n}}\Pi_{h|Y_{1}^{n}} (C_{n}^{c})(1 - \phi_{n}) d\Pi^{U}_{n}(h)
\nonumber
\\&
= \frac{1}{\Pi_{n}(U)}\int 1_{U}(h)\pi(h) \PP_{\theta + h/\sqrt{n}}\Pi_{h|Y_{1}^{n}} (C_{n}^{c})(1 - \phi_{n}) dh
\nonumber
\\&
= \frac{\Pi_{n}(C_{n}^{c})}{\Pi_{n}(U)}\int \PP_{\theta + h/\sqrt{n}}\Pi_{h|Y_{1}^{n}} (U)(1 - \phi_{n}) d\Pi^{C_{n}^{c}}_{n}(h)
\nonumber
\\& \leq \frac{1}{\Pi_{n}(U)} \int_{C_{n}^{c}}\PP_{\theta + h/\sqrt{n}}(1 - \phi_{n})d\Pi_{n}(h)
\nonumber
\\&
\leq
\frac{1}{\Pi_{n}(U)} \int_{\|h\|\geq r_{n}} e^{-Kr_{n}^{2}}d\Pi_{n}(h),
\end{align}
where the last inequality follows from \eqref{eq:bvm:test} and $\frac{1}{\Pi_{n}(U)}$ is bounded above by order $n^{d/2}$ due to the continuity of the density $\pi$ at $\theta_0$; therefore, \eqref{eq:bvm:2} is bounded by
$\int_{\|h\|\geq r_{n}} \exp(-Kr_{n}^{2})dh$, which converges to zero as $n, r_{n} \rightarrow \infty.$

Let $N^{C}(\mu,\Sigma)$ denote the restriction of normal distribution with mean $\mu$ and covariance matrix $\Sigma$ on a set $C$. In the second step of the proof, we want to show that for a ball $C$ with fixed radius $M$,
$$\mathbb{E}_{n,C} \| \Pi^{C}_{h|Y_{1}^{n}} - N^{C}(\Delta_{n,0},I_{\theta_0}^{-1}) \|_{\text{TV}} \rightarrow 0,$$
where $\Delta_{n,0} = \frac{1}{\sqrt{n}} I_{\theta_0}^{-1} l'(\theta_0).$
We have that
\begin{align}
\label{eq:fix1}
&\mathbb{E}_{n,C} \| \Pi^{C}_{h|Y_{1}^{n}} - N^{C}(\Delta_{n,0},I_{\theta_0}^{-1}) \|_{\text{TV}}
\nonumber
\\
&
= 2 \int \int \left\{1 - \frac{dN^{C}(\Delta_{n,0},I_{\theta_0}^{-1})(h)}{1_{C}(h) d\PP_{\theta + h/\sqrt{n}}(Y_{1}^{n})\pi_{n}(h)/\int_{C} d\PP_{\theta + g/\sqrt{n}}(Y_{1}^{n})\pi_{n}(g)dg} \right\}^{+} d\Pi_{h|Y_{1}^{n}}^{C}(dh)\PP_{n,C}(dY_{1}^{n})
\nonumber\\
&
\overset{(i)}{\leq}
2 \int \left\{1 - \frac{1_{C}(h) d\PP_{\theta + g/\sqrt{n}}(Y_{1}^{n})\pi_{n}(g)dN^{C}(\Delta_{n,0},I_{\theta_0}^{-1})(h)}{1_{C}(h) d\PP_{\theta + h/\sqrt{n}}(Y_{1}^{n})\pi_{n}(h)N^{C}(\Delta_{n,0},I_{\theta_0}^{-1})(g)} \right\}^{+}
dN^{C}(\Delta_{n,0},I_{\theta_0}^{-1})(g) d\Pi_{h|Y_{1}^{n}}^{C}(dh)\PP_{n,C}(dY_{1}^{n}),
\end{align}
where $(i)$ follows from Jensen's inequality of $f(x) = (1-x)^{+}$.
By dominated-convergence theorem and $dN^{C}(\Delta_{n,0},I_{\theta_0}^{-1})(g) \leq C\lambda^{C}(g)$ for some constant $C>0$; therefore, it is suffices to show under measure $\PP_{n,C}(dy_{1}^{n})\Pi_{h|Y_{1}^{n}}^{C}(dh) \lambda^{C}(dg)$,
\begin{align}
\label{conv:p}
\frac{d\PP_{\theta + g/\sqrt{n}}(Y_{1}^{n}) \pi_{n}(g) dN^{C}(\Delta_{n,0},I_{\theta_0}^{-1})(h)}{d\PP_{\theta + h/\sqrt{n}}(Y_{1}^{n})\pi_{n}(h) dN^{C}(\Delta_{n,0},I_{\theta_0}^{-1})(g)} \longrightarrow 1 \;\text{in probability},
\end{align}
where $\lambda^{C}$ is the Lebesgue measure on $C$. From the definition of posterior distribution, we have
$$\PP_{n,C}(dy_{1}^{n})\Pi_{h|Y_{1}^{n}}^{C}(dh) \lambda^{C}(dg) = \Pi_{n}^{C}(dh) \PP_{\theta + h/\sqrt{n}}(dy_{1}^{n})\lambda^{C}(dg).$$
Using local asymptotic normality condition \eqref{eq:lan} and continuity of $\pi$, we know that $\PP_{n,C}(dy_{1}^{n})\Pi_{h|Y_{1}^{n}}^{C}(dh) \lambda^{C}(dg)$ is contiguous to $ \lambda^{C}(dh) P_{\theta_0}(dy_{1}^{n})\lambda^{C}(dg);$
therefore, it is suffices to show that under measure $\lambda^{C}(dh) \PP_{\theta_0}(dy_{1}^{n}) \lambda^{C}(dg) $,
$$\frac{d\PP_{n,g}(Y_{1}^{n}) \pi_{n}(g) dN^{C}(\Delta_{n,0},I_{\theta_0}^{-1})(h)}{d\PP_{n,h}(Y_{1}^{n})\pi_{n}(h) dN^{C}(\Delta_{n,0},I_{\theta_0}^{-1})(g)} \longrightarrow 1   \;\text{in probability},$$
which is true under the local asymptotic normality condition \eqref{eq:lan}.
We have for fixed $C$, $$\mathbb{E}_{\theta_0} \| \Pi^{C}_{h|Y_{1}^{n}} - N^{C}(\Delta_{n,0},I_{\theta_0}^{-1}) \|_{TV} \rightarrow 0.$$ Combining this with step one, we have for $n, r_{n} \rightarrow \infty,$
\begin{align}
&\mathbb{E}_{\theta_0} \|\Pi_{h|Y_{1}^{n}}^{C_{n}} - \Pi_{h|Y_{1}^{n}}\|_{\text{TV}} \rightarrow 0,
\nonumber\\&
\|N^{C_{n}}(\Delta_{n,0},I_{\theta_0}^{-1}) - N(\Delta_{n,0},I_{\theta_0}^{-1})\|_{TV} \rightarrow 0.
\end{align}
By diagonal argument, $$\mathbb{E}_{\theta_0} \| \Pi_{h|Y_{1}^{n}} - N(\Delta_{n,0},I_{\theta_0}^{-1}) \|_{TV} \rightarrow 0,$$
which completes the proof.
\hfill $\square$

\section*{Acknowledgements}
Chunlei Wang and Sanvesh Srivastava are partially supported by grants from the Office of Naval Research (ONR-BAA N000141812741) and the National Science Foundation (DMS-1854667/1854662).

\bibliographystyle{Chicago}
\bibliography{papers}

\end{document}